\documentclass[journal]{IEEEtran}
\usepackage[utf8]{inputenc}
\usepackage[margin=3cm]{geometry}

\usepackage{amsmath}
\usepackage{amssymb}
\usepackage{stmaryrd}

\usepackage{macros}

\usepackage{cite}
\usepackage{hyperref}
\usepackage[nottoc,numbib]{tocbibind}

\usepackage{mathtools}
\usepackage[bbgreekl]{mathbbol}
\DeclareSymbolFontAlphabet{\mathbbm}{bbold}
\DeclareSymbolFontAlphabet{\mathbb}{AMSb}
\usepackage{multirow}

\usepackage{caption}

\usepackage{algorithm}
\usepackage[algo2e, linesnumbered]{algorithm2e}

\usepackage[dvipsnames]{xcolor}

\title{Minimal Sparsity 
	for Second-Order Moment-SOS Relaxations \\ of the AC-OPF Problem}
\author{Adrien Le Franc, Victor Magron, Jean-Bernard Lasserre,
		Manuel Ruiz, Patrick Panciatici
	\thanks{Adrien Le Franc, Victor Magron and Jean-Bernard Lasserre
		are with LAAS CNRS, Toulouse, France
		(e-mail: adlefranc@laas.fr)}
	\thanks{Manuel Ruiz and Patrick Panciatici are with 
		R\'eseau de Transport et d’\'Electricit\'e (RTE), Paris, France}}

\date{\today}

\begin{document}

\maketitle

\begin{abstract}
	
	AC-OPF (Alternative Current Optimal Power Flow)
	aims at minimizing the operating costs of a power grid
	under physical constraints on voltages and power injections.
	Its mathematical formulation results in a nonconvex polynomial optimization
	problem which is hard to solve in general,
	but that can be tackled 
	by a sequence of SDP
	(Semidefinite Programming) relaxations
	corresponding to the steps of 
	the moment-SOS (Sums-Of-Squares) hierarchy.
	Unfortunately, the size of these SDPs 
	grows drastically in the hierarchy,
	so that even second-order relaxations
	exploiting the correlative sparsity pattern of AC-OPF
	are hardly numerically tractable for large
	instances --- with thousands of power buses.
	Our contribution lies in a new sparsity
	framework, termed minimal sparsity, inspired
	from the specific structure of power flow
	equations.
	Despite its heuristic nature, 
	numerical examples show that minimal sparsity allows the computation of
	highly accurate second-order moment-SOS relaxations
	of AC-OPF, while requiring far less computing time 
	and memory resources than the standard correlative sparsity pattern. 
	Thus, we manage to compute second-order relaxations on test cases
	with about 6000 power buses, which we believe  
	to be unprecedented.

\end{abstract}


\section*{Nomenclature}

\begin{itemize}
	
	\item For a finite set $F$,
	we write $\card{F}$ for its cardinality.
	
	\item For a complex number $z \in \CC$,
	we write
	$\angle z$ for its angle;
	$\abs{z}$ for its magnitude;
	$z^*$ for its complex conjugate;
	$\Real{z}$ for its real part;
	and $\Imaginary{z}$ for its imaginary part.
	
	\item For a pair of integers $(a,b) \in \NN^2$
	with $a \leq b$, we write
	$\ic{a,b}$ for the sequence $\na{a, a+1, \ldots, b}$.
	
	\item For an $N\times N$ real symmetric matrix
	$M \in \SYM^N$, $M \succeq 0$ means that $M$
	is positive semidefinite (PSD). 
	
	\item For real matrices 
	$(A,B) \in (\RR^{N\times N})^2$, we write
	$\proscal{A}{B}$ for the Frobenius inner product
	between $A$ and $B$.
	
	\item For a 
	polynomial function $f : \RR^N \to \RR$
	that decomposes as $\sum_{\alpha \in \NN^N} f_\alpha x^\alpha$
	in the standard monomial basis, we denote its support
	by $\supp(f) = \defset{\alpha \in \NN^N}{f_\alpha \neq 0}$,
	and the set of variables involved in $f$ by
	$\var(f) = \defset{n \in \ic{1,N}}{ \exists \alpha \in \supp(f)
		\eqsepv \alpha_n \neq 0 }$.
	
\end{itemize}

Whether $\card{\cdot}$ denotes cardinality or magnitude 
is always clear from context. 

\section{Introduction}

\IEEEPARstart{T}{he} AC-OPF 
(Alternative Current - Optimal Power Flow) problem
plays a central role for the management of AC power grids,
but remains highly challenging to solve.
Indeed, AC-OPF formulates as a nonconvex optimization program,
and real scale instances typically have 
thousands of decision variables \cite{babaeinejadsarookolaee2019power}.
A common way of addressing this problem is to compute
a local solution with a nonlinear solver, but the solution obtained might not be globally optimal \cite{bukhsh2013local}.
Therefore, a vast body of literature have concentrated on
relaxations of the original problem to compute lower
bounds so as to estimate the quality of a local solution.
We refer to \cite{molzahn2019survey} for a recent survey
of such relaxations.

In this paper,
we follow the approach of \cite{josz2014application, molzahn2014moment}
to compute lower bounds of AC-OPF instances
based on the moment-SOS (Sums-Of-Squares) hierarchy \cite{lasserre2001global}.
In this framework, we consider a sequence of SDP
(Semidefinite Programming) relaxations
of the AC-OPF problem, whose values
converge monotonously to the optimal value of the original problem.
This sequence generically converges in a finite
number of steps \cite{nie2014optimality} and, experimentally,
the second-order relaxation already achieves this convergence for 
most of the AC-OPF test cases reported in \cite{josz2014application, molzahn2014moment, gopinath2020proving}.
However, the size of the SDPs involved in the hierarchy grows drastically
with the number of AC-OPF variables and with the order of the relaxation,
so that this method becomes rapidly challenging from a numerical perspective.

Exploiting the correlative sparsity pattern \cite{waki2006sums} 
of AC-OPF  has led to a hierarchy of SDP 
relaxations of reduced size, and has helped a lot to scale to larger networks \cite{molzahn2014sparsity}. 
See also the recent survey \cite{sparsebook} for several applications of sparse polynomial optimization. 
Yet, even these sparse second-order relaxations yield out-of-memory errors 
for some instances 
with about a hundred of power buses
on a computer with 125 GB of RAM \cite{gopinath2020proving}.
Regarding instances with thousands of buses, the most scalable
approach to our knowledge seems to be the recent
correlative-term sparsity framework \cite{wang2022cs},
which enables the computation of partial sparse second-order moment relaxations
for large scale power grids \cite{wang2021certifying}. 
Nevertheless,
the accuracy of second-order relaxations
remains attractive, as it provides lower bounds
certifying 0.00\% of optimality gaps on almost all
tractable instances of~\cite{gopinath2020proving}.

{\bf Contribution:}
we focus our attention on sparse second-order moment-SOS relaxations
of AC-OPF based on correlative sparsity.
In the standard case, the size of the matrix variables involved 
is ruled by a family of subsets of optimization variables 
obtained by an algorithmic routine --- computing the maximal cliques of a chordal graph.
This approach enforces the RIP
(Running Intersection Property),
which ensures the convergence of the sparse hierarchy. 
We refer to \cite{lasserre2006convergent} for technical details,
including the formal definition of the RIP. 

%

Our contribution lies in a new sparsity framework, 
that we call \textit{minimal sparsity}.
%
%
Compared  with  the  standard  approach,  we chose  to  relax  the  RIP  so  as  to  gain  control  on the  size  of  matrix  variables.
Our definition of minimal sparsity is inspired by the specific
structure of the power flow equations and results in sparse second-order
relaxations that have smaller matrix variables --- which is generally preferred 
by SDP solvers based on interior-point methods \cite{nesterov1994interior}.
Therefore, our approach is heuristic
--- as we relax the RIP, convergence to the global minimum is not guaranteed any more ---
although we may still
easily extract a global optimal solution of AC-OPF
if the moment matrices of our relaxations are rank-one.
In spite of this heuristic nature,
we report that
minimal sparsity yields highly accurate lower bounds on practical test cases.
Moreover, the second-order relaxations obtained scale much better than their clique-based
counterparts, allowing us to handle AC-OPF instances with thousands of buses. 

The paper is organized as follows.
First, in~\S\ref{sec:sparse_moment_SOS_in_ACOPF},
we recall background notions on sparse moment-SOS hierarchies 
and their application to AC-OPF.
Second, in~\S\ref{sec:minimal_sparsity_for_scalable_ACOPF},
we introduce our new minimal sparsity framework.
Third, in~\S\ref{sec:numerical_examples}, we illustrate the strengths of 
minimal sparsity by computing second-order moment-SOS relaxations
of the AC-OPF problem on various numerical test cases.

\section{Sparse moment-SOS relaxations for the AC-OPF problem}
\label{sec:sparse_moment_SOS_in_ACOPF}%

First, in~\S\ref{sec:AC-OPF}, we recall the formulation of the AC-OPF problem.
Second, in~\S\ref{sec:sparse_moment_hierarchies},
we review basic concepts of moment-SOS hierarchies
and their applications 
to certify global optimality in AC-OPF.
Third, in~\S\ref{sec:sparse_relaxations},
we present background notions on sparse moment-SOS relaxations.

\subsection{The AC-OPF problem}
\label{sec:AC-OPF}

\begin{figure}[htbp]
	\centering
	\includegraphics[width=0.4\textwidth]{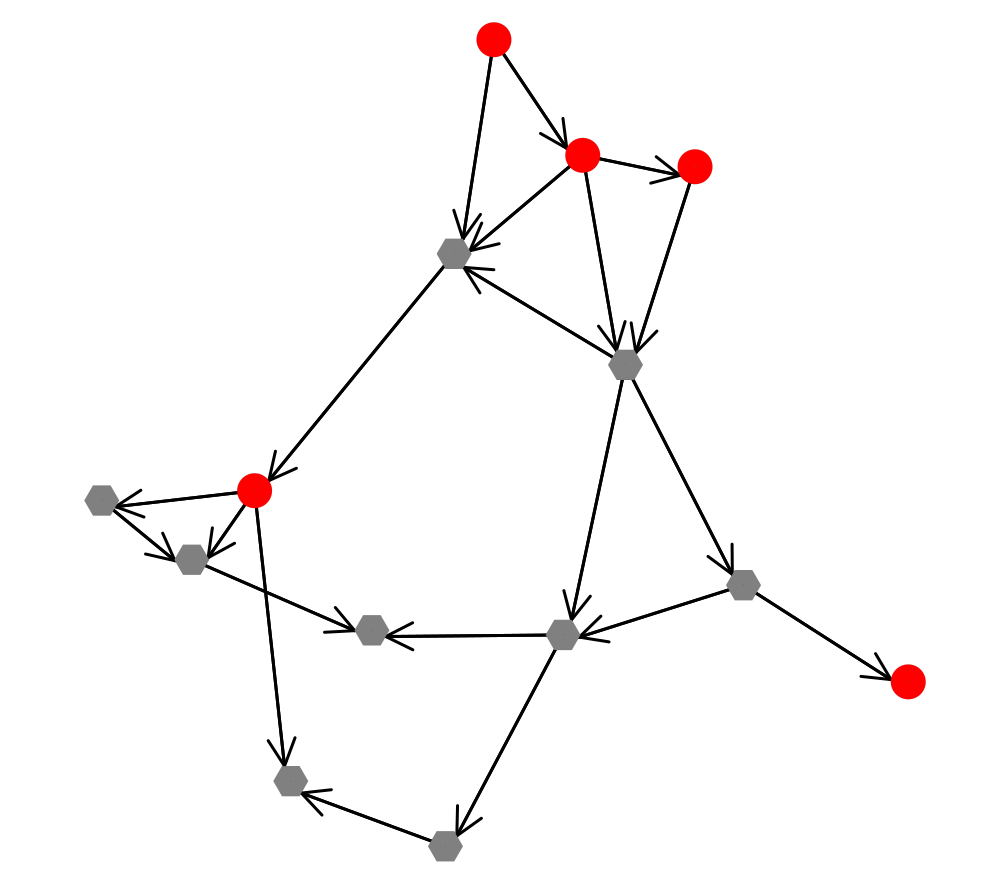}
	\caption{Example of $(\Nodes, \Edges)$ graph model
		for PGLib's case 14 IEEE:
		nodes represent buses and edges
		represent power lines with their conventional
		power flow orientation. 
		Red circle node markers
		highlight buses with power generators.}
	\label{fig:grid}
\end{figure}

In the AC-OPF problem, we aim at minimizing
the operating costs of a power grid while satisfying
power flow balance equations and infrastructure constraints.  
We model the grid by a directed graph $(\Nodes, \Edges)$
where nodes $\Nodes$ represent buses and edges $\Edges$
represent power lines. 
Line orientations model the asymmetry of power flow
along transmission lines in AC power grids.
A subset $\Generators \subseteq \Nodes$ of nodes 
highlights buses with generating power units.
An illustrative $(\Nodes, \Edges)$ example based on PGLib's case 14 IEEE
is given in Figure~\ref{fig:grid}.
For the sake of clarity, we assume here that at most one single line
can connect two buses $(i,j) \in \Nodes^2$.
Parallel lines can be modeled by adequate edge labeling
as in~\cite[Model 1]{babaeinejadsarookolaee2019power}.

Formally, AC-OPF amounts to solving the following optimization problem:

\begin{subequations}
	\label{eq:AC-OPF}
	\begin{align}
	&\min_{\substack{\voltage \in \CC^{\card{\Nodes}} \\
			\power \in \CC^{\card{\Generators}} \\
			\powerflow \in \CC^{2\card{\Edges} } }} \quad
	 \sum_{g \in \Generators} C_{2,g} \Real{\power_g}^2 
	+ C_{1,g}\Real{\power_g} + C_{0,g} \eqfinv
	\label{eq:power_gen_costs}
	\intertext{s.t.}
	&\forall i \in \RefNodes \eqsepv \angle\voltage_i = 0 \eqfinpv
	\label{eq:reference_voltage}
	\\
	 &\forall g \in \Generators \eqsepv \minpower_g \leq \power_g \leq \maxpower_g \eqfinpv
	\label{eq:generator_limits}
	\\
	&\forall i \in \Nodes \eqfinv \notag \\
	&\minmodulus_i \leq \abs{\voltage_i} \leq \maxmodulus_i 
	\eqfinv 
	\label{eq:voltage_limits}
	\\
	 &\sum_{g \in \Generators(i)} \power_g - \load_i - \bp{\shuntadmittance_i}^*
	\abs{\voltage_i}^2
	= \sum_{j \in \Neighbors{i}} \powerflow_{i,j}
	\eqfinpv
	\label{eq:power_flow}%
	\\
	&\forall (i,j) \in \Edges \eqfinv \notag
	\\
	 &\powerflow_{i,j} = \bp{\admittance_{i,j} + \admittance_{i,j}^c}^*
	\frac{\abs{\voltage_i}^2}{\abs{\ratio_{i,j}}^2}
	- \admittance_{i,j}^*\frac{\voltage_i\voltage_j^*}{\ratio_{i,j}}
	\eqfinv 
	\label{eq:line_power_1}%
	\\
	 &\powerflow_{j,i} = \bp{\admittance_{i,j} + \admittance_{j,i}^c}^*
	\abs{\voltage_j}^2
	- \admittance_{i,j}^*\frac{\voltage_i^*\voltage_j}{\ratio_{i,j}^*}
	\eqfinv 
	\label{eq:line_power_2}%
	\\
	 &\abs{\powerflow_{i,j}} \leq \maxpowerflow_{i,j}
	\eqsepv 
	\abs{\powerflow_{j,i}} \leq \maxpowerflow_{i,j}
	\eqfinv \label{eq:line_constraints}
	\\
	 &\minangle_{i,j} \leq \angle (\voltage_i\voltage_j^*) \leq \maxangle_{i,j}
	\eqfinp
	\label{eq:angle_constraints}
	\end{align}
\end{subequations}

In the above formulation, 
lower case letters are used for decision variables
and capital letters refer to constant parameters.
The original decision variables
are the bus voltages $\voltage = \na{\voltage_i}_{i \in \Nodes}$
and the power generation values $\power = \na{\power_g}_{g \in \Generators}$.
Additionally, for every edge $(i, j) \in \Edges$,
we introduce $\powerflow_{i,j}$
for the power flow from bus $i$ to bus $j$
and $\powerflow_{j,i}$ for the power flow from bus $j$ to bus $i$.

We now provide physical interpretations 
for the objective and constraints of Problem~\eqref{eq:AC-OPF}:

\begin{itemize}
	
	\item We minimize power generation costs~\eqref{eq:power_gen_costs},
	which are assumed to only depend on the real part of $\power_g$,
	for $g \in \Generators$ --- that is, on active power generation ---
	with parameters $\np{C_{0,g}, C_{1,g}, C_{2,g}} \in \RR^3$.
	
	\item In constraint~\eqref{eq:reference_voltage},
	we set the voltage angle of some reference
	buses $\RefNodes \subseteq \Nodes$ to zero
	to address the rotational invariance of voltage solutions.
	
	\item In constraints~\eqref{eq:generator_limits}-\eqref{eq:voltage_limits},
	we impose bounds $\np{\minpower_g, \maxpower_g} \in \CC^2$ on the real and imaginary parts 
	of the generated power $\power_g$, for $g \in \Generators$,
	and bounds 
	$\np{\minmodulus_i, \maxmodulus_i} \in \RR^2_+$
	on the magnitude of the bus voltage 
	$\voltage_i$, for $i \in \Nodes$.
	
	\item In constraint~\eqref{eq:power_flow},
	we enforce the balance of power flows at
	every bus $i \in \Nodes$.
	The balance equation involves
	power generations $\power_g$ for $g$
	in the (possibly empty)
	set $\Generators(i) \subseteq \Generators$
	of generators at bus $i$;
	power flows $\powerflow_{i,j}$
	for $j$ in the set $\Neighbors{i} \subseteq \Nodes$
	of neighbors of bus $i$;
	the load $\load_i \in \CC$ 
	and a shunt admittance term with $\shuntadmittance_i \in \CC$.
	
	\item In constraints~\eqref{eq:line_power_1}-\eqref{eq:line_power_2},
	we give the expression of
	power flows ($\powerflow_{i,j}, \powerflow_{j,i}$)
	along every line $(i,j) \in \Edges$,
	following the $\Pi$-circuit branch model
	with parameters 
	$(\admittance_{i,j}, \admittance_{i,j}^c,
	\admittance_{j,i}^c, \ratio_{i,j}) \in \CC^4$
	detailed in~\cite[Appendix B]{babaeinejadsarookolaee2019power}.
	
	\item In constraints~\eqref{eq:line_constraints}-\eqref{eq:angle_constraints},
	we impose a thermal limit $\maxpowerflow_{i,j} \in \RR_+$
	on power flows and voltage angle difference bounds
	$(\minangle_{i,j}, \maxangle_{i,j}) \in \RR^2$ for every line $(i,j) \in \Edges$.
		
\end{itemize}

Due to nonlinear equality and nonconvex inequality
constraints, Problem~\eqref{eq:AC-OPF} is nonconvex,
and hard to solve in general~\cite{bienstock2019strong}.

\subsection{SDP lower bounds via moment-SOS hierarchies}
\label{sec:sparse_moment_hierarchies}

Following the approach of \cite{josz2014application, molzahn2014moment}, AC-OPF can be cast
as a POP (Polynomial Optimization Problem)
to benefit from powerful results of the moment-SOS
hierarchy. 
We introduce notations for such a reformulation
of Problem~\eqref{eq:AC-OPF} and recall some fundamental
properties of moment-SOS relaxations.

\paragraph{From AC-OPF to POP}

by considering separately the real and imaginary parts
of voltage and power generation
variables of Problem~\eqref{eq:AC-OPF},
we obtain $\nmax = 2(\card{\Nodes} + \card{\Generators})$
real variables $x \in \RR^\nmax$
(power flow variables are omitted by injecting 
\eqref{eq:line_power_1}-\eqref{eq:line_power_2} into~\eqref{eq:power_flow}).
The correspondence between AC-OPF and POP variables
is formalized by two bijective mappings
\begin{subequations}
\begin{equation}
	\bijR : \Generators \cup \Nodes \to \ic{1, \frac{\nmax}{2}}
	\eqsepv 
	\bijIm : \Generators \cup \Nodes \to \ic{\frac{\nmax}{2}+1, \nmax}
	\eqfinv
\end{equation}
so that
\begin{align}
	x_{\bijR(g)} &= \Real{\power_g} \eqsepv 
	x_{\bijIm(g)} = \Imaginary{\power_g}
	\eqsepv
	\forall g \in \Generators \eqfinv
	\\
	x_{\bijR(i)} &= \Real{\voltage_i} \eqsepv 
	x_{\bijIm(i)} = \Imaginary{\voltage_i}
	\eqsepv
	\forall i \in \Nodes \eqfinp
\end{align}
\end{subequations}

Then, we observe that every constraint in~\eqref{eq:reference_voltage}-\eqref{eq:angle_constraints}
can be equivalently formulated as an equality or inequality
constraint defined with a multivariate polynomial 
in $x \in \RR^\nmax$. To perform this reformulation,
constraints~\eqref{eq:voltage_limits}
and~\eqref{eq:line_constraints} need to be squared
and constraint~\eqref{eq:angle_constraints}
needs to be transformed as detailed e.g. 
in~\cite[\S5.1.2]{bienstock2020mathematical}.
Thus, by introducing $\kmax+1$ appropriate  
real multivariate polynomial functions
$\na{f_k}_{k \in \ic{0,\kmax}}$,
Problem~\eqref{eq:AC-OPF} can be written as a POP:
\begin{subequations}
\label{eq:POP}
\begin{align}
\rho &= \min_{x \in \ConstraintSet} \f_0(x) \eqsepv \text{ where } \enspace \\  
\ConstraintSet &= \defset{x \in \RR^\nmax}{\f_k(x) \geq 0 
	\eqsepv \forall k \in \ic{1,\kmax}} \eqfinp
\end{align}
\end{subequations}

\paragraph{The Moment-SOS hierarchy}

despite its potential nonconvexity, the optimal value
of Problem~\eqref{eq:POP} can be approximated
--- and often exactly computed --- by the moment-SOS hierarchy \cite{lasserre2001global}.
In this framework, we consider two 
sequences of SDPs,
starting from a minimal order 
$r_0 = \max \na{d_k}_{k \in \ic{0,K}}$
where $d_k = \big\lceil \frac{\text{deg}(f_k)}{2} \big\rceil$.
The moment hierarchy is defined by a sequence of SDPs indexed by $r \in \ic{r_0, +\infty}$:
\begin{subequations}
	\begin{align}
	\val_r = & \min_{y} \quad  
	\sum_{\alpha \in \supp(f_0)} f_{0,\alpha} y_\alpha 
	\eqfinv
	\\
	\text{s.t.} \quad 
	& \MomentMatrix_r(y) \succeq 0 \eqfinv
	\label{eq:moment_matrix}%
	\\
	&  \MomentMatrix_{r-d_k}(f_k y) \succeq 0 \eqsepv
	\forall k \in \ic{1, \kmax} \eqfinv
	\label{eq:localization_matrices}%
	\\
	& y_{\zero} = 1 \eqfinp
	\end{align}
	\label{eq:moment_hierarchy}%
\end{subequations}
The entries of the so-called \textit{pseudo-moment} variable vector $y$ in Problem~\eqref{eq:moment_hierarchy}
are indexed by elements of
the truncated monomial basis 
$\na{x^\alpha}_{\alpha \in \NN_{2r}^\nmax}$,
where 
$\NN_{r}^\nmax = \ba{\alpha \in \NN^\nmax 
\ | \ \sum_{n \in \ic{1:\nmax}} \alpha_n \leq r}$ for $r \in \NN$.
Subsequently, the moment matrix in~\eqref{eq:moment_matrix} 
and the localization matrices
in~\eqref{eq:localization_matrices}
are expressed as
\begin{subequations}
	\label{eq:moment_and_localizing_matrices_definition}%
	\begin{align}
	\MomentMatrix_r(y) &= \np{y_{\alpha+\beta}}_{\alpha, \beta \in \NN_r^\nmax}
	\eqfinv
	\\
	\MomentMatrix_{r-d_k}(f_k y) &= 
	\Bp{\sum_{\gamma \in \supp(f_k) }f_{k,\gamma} y_{\alpha+\beta+\gamma}}_{\alpha, \beta \in \NN_{r-d_k}^\nmax}
	\eqfinp
	\end{align}
\end{subequations}
These matrices have entries that are linear in the ones of $y$, 
so that we can write
$\MomentMatrix_r(y) = \sum_{\alpha \in \NN_{2r}^\nmax} 
\Mat_{0,\alpha} y_\alpha$
and
$\MomentMatrix_{r-d_k}(f_k y) = \sum_{\alpha \in \NN_{2r}^\nmax} 
\Mat_{k,\alpha} y_\alpha$
by introducing adequate matrices 
$\na{\Mat_{k,\alpha}}_{\alpha \in \NN_{2r}^\nmax}$
for all $k \in \ic{0, \kmax}$.

By considering the dual of~\eqref{eq:moment_hierarchy}, we obtain the SOS hierarchy of SDPs indexed by $r \in \ic{r_0, +\infty}$:
\begin{subequations}
	\label{eq:sos_hierarchy}
	\begin{align}
	\valSOS_r =  
	& \max_{\G, t} 
	\quad t \eqfinv
	\\ 
	\text{s.t.} \quad
	& f_{0, \zero} - t = \sum_{k\in\ic{0,\nmax}}
	\proscal{\Mat_{k,\zero}}{\G_k} 
	\eqfinv 
	\label{eq:SOS_decomposition_0}
	\\
	& f_{0, \alpha}  = \sum_{k\in\ic{0,\nmax}}
	\proscal{\Mat_{k,\alpha}}{\G_k} 
	\eqsepv
	\forall \alpha \in \NN_{2r}^\nmax \setminus
	\na{\zero} \eqfinv 
	\label{eq:SOS_decomposition}
	\\
	& \G_k \succeq 0 \eqsepv \forall k \in \ic{0, K} \eqfinp 
	\end{align}
\end{subequations}
In the context of AC-OPF, 
strong duality holds between 
Problems~\eqref{eq:moment_hierarchy} and \eqref{eq:sos_hierarchy} 
(see \cite{josz2014application}), 
and the nondecreasing sequences of lower bounds
$\na{\val_r}_{r \geq r_0}$ and $\na{\valSOS_r}_{r \geq r_0}$
converge to the value $\val$ of the POP~\eqref{eq:POP}
(see \cite{lasserre2001global}).
However, the sizes of the corresponding SDP relaxations grow drastically 
with the values of $\nmax$ and $r$, as 
the largest Gram matrix $\G_0$ in~\eqref{eq:sos_hierarchy}
and the moment matrix $\MomentMatrix_r(y)$ in~\eqref{eq:moment_hierarchy}
are of size $\abs{\NN_{r}^\nmax} = \binom{\nmax+r}{r}$.

\subsection{Sparse relaxations}
\label{sec:sparse_relaxations}%



One way to bypass 
the curse of dimensionality
mentioned hereabove 
is to exploit the sparsity
of AC-OPF, as initially suggested in \cite{molzahn2014sparsity}.
In the context of the moment hierarchy, sparsity
consists in reducing the dimension
of the search space
of Problem~\eqref{eq:moment_hierarchy}
by selecting a subset of monomials
in $\na{x^\alpha}_{\alpha \in \NN_{2r}^\nmax}$
for indexing the pseudo-moment variable vector $y$.
We concentrate on the correlative sparsity
pattern~\cite{waki2006sums}
which introduces a hierarchy of sparse moment relaxations:
%
\begin{subequations}
	\label{eq:sparse_moment_hierarchy}%
	\begin{align}
	\val_r(\Indices) =&  \min_{ y } 
	\quad  
	\sum_{\alpha \in \supp(f_0)} f_{0,\alpha} y_\alpha \eqfinv
	\\
	\text{s.t.} \quad & \MomentMatrix_{r}(y ; \Indices_p) \succeq 0 \eqsepv
	\forall p \in \ic{1,\pmax} \eqfinv \label{eq:sparse_moment_matrix}
	\\
	& \MomentMatrix_{r-d_k}(f_k y ; \Indices_p) \succeq 0 \eqsepv
	\forall k \in \Constraints_p \eqsepv 
	\label{eq:sparse_localization_matrices}
	\\
	& \vspace{3cm} \forall p \in \ic{1,\pmax} \eqfinv  \notag
	\\
	&  y_{\zero} = 1 \eqfinp
	\end{align}
\end{subequations}
Problem~\eqref{eq:sparse_moment_hierarchy}
is parameterized by
a family of subsets of $\ic{1,\nmax}$,
denoted
$\Indices = \na{\Indices_p}_{\in \ic{1,\pmax}}$,
and satisfying
$\cup_{p \in \ic{1,\pmax}} \Indices_p = \ic{1,\nmax}$.
The constraints $\na{f_k}_{k \in \ic{1,\kmax}}$ 
are distributed over a partition $\na{\Constraints_p}_{\in \ic{1,\pmax}}$ of $\ic{1,\kmax}$ such that
for all $p \in \ic{1,\pmax}$ and $k \in \Constraints_p$,
$\var(f_k) \subseteq \Indices_p$.
Then,
for $p \in \ic{1,\pmax}$,
the sparse moment and localization matrices in~\eqref{eq:sparse_moment_matrix}-\eqref{eq:sparse_localization_matrices}
are defined after~\eqref{eq:moment_and_localizing_matrices_definition}
by selecting only rows and columns indexed by monomials
in $\na{x^\alpha}_{\alpha \in \NN_{2r}^\nmax}$
satisfying $\var(x^\alpha) \subseteq \Indices_p$.
Naturally, the dual of Problem~\eqref{eq:sparse_moment_hierarchy}
gives rise to a sparse SOS hierarchy,
whose sequence of bounds is introduced as
$\na{\valSOS_r(\Indices)}_{r \geq r_0}$.

We remind that the choice of the subsets $\Indices$ 
is of paramount importance.
On the practical side, the cardinalities
of these subsets control
the sizes of the matrices in~\eqref{eq:sparse_moment_matrix}-\eqref{eq:sparse_localization_matrices}.
In general, 
the smaller these matrices, 
the better the numerical performances of SDP solvers,
especially for those based on interior-point methods \cite{nesterov1994interior}.
On the theoretical side, the bounds 
$\na{\val_r(\Indices)}_{r \geq r_0}$
are not guaranteed to converge to 
the value $\val$ of the POP~\eqref{eq:POP} 
for any choice of $\Indices$.
The most favorable case is when the subsets $\Indices$
satisfy the RIP (Running Intersection Property)
where asymptotic convergence
is preserved~\cite{lasserre2006convergent}. 
These considerations on the design of $\Indices$
are further investigated
in the next section.

\section{Minimal sparsity for scalable AC-OPF relaxations}
\label{sec:minimal_sparsity_for_scalable_ACOPF}%

We recall basic notions of clique-based sparsity
and expose some of its limitations 
regarding computing scalability 
in~ \S\ref{sec:clique_based_sparsity}.
As an alternative, we introduce
our minimal sparsity pattern in~\S\ref{sec:minimal_sparsity}.
We further detail a method to control
the cardinalities of the subsets $\Indices$
in~\S\ref{sec:finer_control}.

\subsection{Clique-based sparsity and its limitations}
\label{sec:clique_based_sparsity}

We recall how to compute clique-based subsets $\Indices$
and discuss some limitations of this approach.

\paragraph{Clique-based subsets}

the design of subsets $\Indices$
satisfying the RIP is usually based on the following algorithmic
routine.

\begin{itemize}
	\item[$(i)$] 
	First, we define the \textit{correlative sparsity pattern graph}
	$\np{\NodesCSP, \EdgesCSP}$.
	In this graph,
	nodes $\NodesCSP$ represent the $\nmax$ variables
	of the POP~\eqref{eq:POP}
	and undirected edges $\EdgesCSP$
	account for products between variables
	in the polynomial functions $\na{f_k}_{k \in \ic{1,\kmax}}$: 
	an edge $(n_1, n_2) \in \EdgesCSP$ indicates that
	there exists $\alpha \in \cup_{k \in \ic{0,\kmax}} \supp(f_k)$
	such that $\na{n_1, n_2} \subseteq \var(x^\alpha)$.
	
	\item[$(ii)$]
	Second, we perform a chordal extension of
	$\np{\NodesCSP, \EdgesCSP}$. We recall that a graph
	is chordal if each of its cycle of length four ou greater
	has a chord. Therefore, chordal extension adds new edges,
	resulting in a new graph $\np{\NodesCSP, \EdgesCSPChordal}$,
	where $\EdgesCSP \subseteq \EdgesCSPChordal$.
	
	\item[$(iii)$]	
	Third, we define subsets $\Indices^\text{c}$
	as the maximal cliques of the chordal graph
	$\np{\NodesCSP, \EdgesCSPChordal}$.
	We recall that a clique is a complete subgraph
	of $\np{\NodesCSP, \EdgesCSPChordal}$, and that it is maximal when it cannot be augmented
	by adding an adjacent node.
			
\end{itemize}

The clique-based subsets $\Indices^\text{c}$
satisfy the RIP, and thus
ensure the convergence of the correlative
sparse moment-SOS hierarchy~\cite{lasserre2006convergent}.
%

\paragraph{Limitations of clique-based sparsity}

the above routine for designing the subsets $\Indices$
gives a systematic way to reduce the computing
burden of the dense relaxation~\eqref{eq:moment_hierarchy}.
However, for large AC-OPF instances, even the sparse
relaxation~\eqref{eq:sparse_moment_hierarchy}
can be numerically challenging.
Experimentally, \cite{gopinath2020proving}
report that the second-order sparse moment relaxation
triggers an out-of-memory error 
on a computer allowed with 125 GB of RAM
for PGLib's instances 89 Pegase and 162 IEEE.

Therefore, some works concentrate
on improving the algorithmic routine $(i)-(iii)$
to reduce memory usage and computing time for solving~\eqref{eq:sparse_moment_hierarchy}.
In particular, \cite{molzahn2013implementation, sliwak2020clique}
propose clique merging strategies as a post-processing of $(iii)$.
This line of work has allowed up to $\times$3 decreases
in solving time for first-order relaxations \cite{molzahn2013implementation}.
However, extensions to second-order relaxations
seem much less effective \cite[\S 4.5]{sliwak2021resolution}.
We believe that it is due to the iteration complexity
of interior-point SDP solvers,
which typically perform operations that scale cubically 
with the size of the largest SDP matrix \cite{nesterov1994interior}.
We recall that the largest
matrix in~\eqref{eq:sparse_moment_hierarchy} is of size
\begin{equation}
	\binom{m +r}{r}
	\eqsepv
	\text{ where } m = \max_{p \in \ic{1,\pmax}} \card{\Indices_p}
	\eqfinv 
\end{equation}
hence the importance of moderating the cardinalities
of the subsets in $\Indices$ to alleviate memory requirements and 
computing time in second-order sparse relaxations.
Clearly, a clique merging strategy is not meant to
reduce these cardinalities,
and therefore does not address what 
we identify as the principal bottleneck in second-order
sparse relaxations.

\subsection{Minimal sparsity}
\label{sec:minimal_sparsity}

We introduce minimal subsets $\IndicesMS$
to address the principal 
limitations faced with clique-based sparsity
in AC-OPF.
Our definition builds on the specific structure of
power flow equations:
for each bus $i \in \Nodes$,
we select the minimal group of POP variables
required to write the power flow balance equation
at bus $i$.
This results in $\pmax = \card{\Nodes}$ subsets
given by
\begin{subequations}
\begin{align}
\IndicesMS_{\#i} =& \na{\bijR(i), \bijIm(i)}
\bigcup_{j \in \Neighbors{i}} \na{\bijR(j), \bijIm(j)} \notag \\
&\bigcup_{g \in \Generators(i)} \na{\bijR(g), \bijIm(g)}
\eqfinv \label{eq:minimal_sparsity}
\end{align}
where, assuming an arbitrary order on buses $\Nodes$,
we denote by $\#i \in \ic{1,\pmax}$ the position of 
bus $i \in \Nodes$.
In term of correspondence between POP and AC-OPF 
formulations, we obtain the following relationship:
\begin{align}
	\na{x_n}_{n \in \IndicesMS_{\#i}}
	=& \na{\Real{\voltage_i}, \Imaginary{\voltage_i}}
	\bigcup_{j \in \Neighbors{i}} 
	\na{\Real{\voltage_j}, \Imaginary{\voltage_j}} \notag \\
	&\bigcup_{g \in \Generators(i)} 
	\na{\Real{\power_g}, \Imaginary{\power_g}}
	\eqfinp
	\label{eq:minimal_sparsity_variables}%
\end{align}
\end{subequations}
The above expression highlights 
that in $\IndicesMS_{\#i}$, we select 
the minimal amount of AC-OPF variables
that we need to write constraints~\eqref{eq:power_flow}
and~\eqref{eq:line_power_1}-\eqref{eq:line_power_2}
at bus $i \in \Nodes$ for Problem~\eqref{eq:AC-OPF}.

\begin{figure}[htbp]
	\centering
	\includegraphics[width=0.5\textwidth]{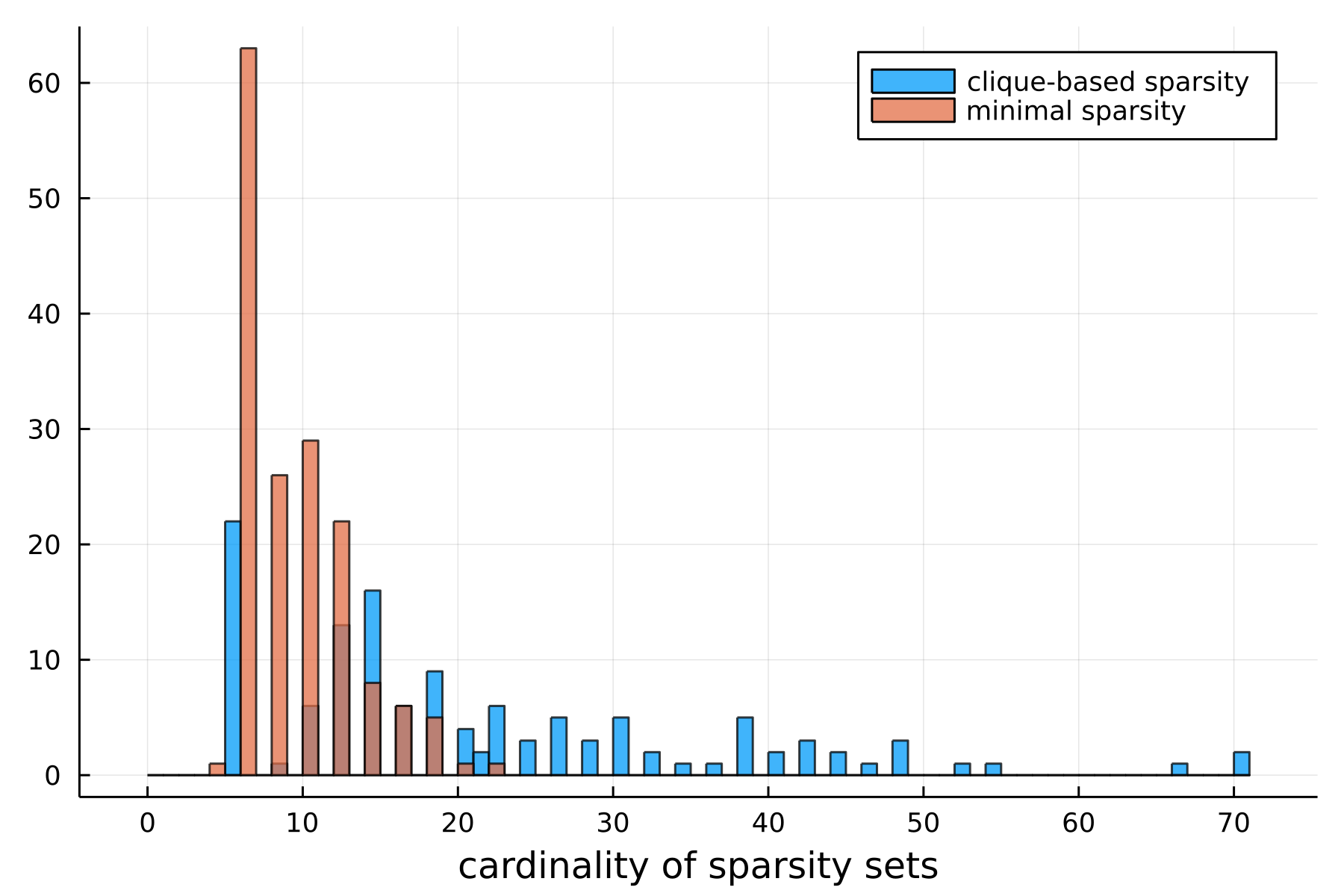}
	\caption{Histogram of the cardinalities
		of clique-based subsets $\IndicesCB$
		(blue color) and minimal subsets $\IndicesMS$
		(red color)
		for PGLib's case 162 IEEE}
	\label{fig:sparsity_sets}
\end{figure}

Minimal sparsity entails a trade-off
between the number of subsets in $\Indices$ 
and their cardinalities. 
We illustrate this trade-off 
by comparing clique-based
subsets $\IndicesCB$ and minimal subsets $\IndicesMS$ 
for PGLib's case 162 IEEE.
We compute the chordal extension 
$\np{\NodesCSP, \EdgesCSPChordal}$ and its maximal cliques
using the \textit{greedy fillin} heuristic 
implemented in the TSSOS package \cite{magron2021tssos},
as this heuristic is expected to yield
smaller average clique numbers
than other standard heuristics \cite{bodlaender2010treewidth}. 
The histogram of the cardinalities of sets
for both sparsity patterns is given in Figure~\ref{fig:sparsity_sets}.
For case 162 IEEE,
$\IndicesCB$ has $\pmax = 126$ sets,
the largest of which has 70 variables,
whereas $\IndicesMS$ has $\pmax = 162$ sets
with at most 22 variables.
Consequently, the sparse moment-SOS relaxations
written with $\IndicesMS$ have a larger amount
of PSD constraints-matrices than clique-based
relaxations, but their dimensions
are much smaller.
In general, this situation is preferred by SDP solvers based 
on interior-point methods \cite{nesterov1994interior}. 

\subsection{Finer control on the size of subsets}
\label{sec:finer_control}%

If the graph $(\Nodes, \Edges)$ has nodes
with a high number of neighbors,
the minimal subsets $\IndicesMS$
defined by~\eqref{eq:minimal_sparsity}
may still have large cardinalities.
Assuming that we wish to impose a maximal cardinality threshold
$\Imax$ for the subsets $\Indices^\text{m}$,
we propose a modification of the AC-OPF Problem~\eqref{eq:AC-OPF}
and of the minimal subsets $\Indices^\text{m}$
to meet this requirement.

In our approach, when $\card{\IndicesMS_{\#i}} > \Imax$
at some bus $i \in \Nodes$,
we split neighboring buses $\Neighbors{i}$
into a partition
$\na{\NeighborsA{i}}_{a \in \Additional(i)}$,
where the set $\Additional(i)$
is introduced to index additional
complex variables $\na{\additional_{i,a}}_{a \in \Additional(i)}$
for the AC-OPF Problem~\eqref{eq:AC-OPF}.
Then, we rewrite the power flow
equation~\eqref{eq:power_flow} at bus $i$ as
\begin{subequations}
\begin{align}
	&\sum_{g \in \Generators(i)} \power_g - \load_i - \bp{\shuntadmittance_i}^*
	\abs{\voltage_i}^2
	= \sum_{a \in \Additional(i)} \additional_{i,a}
	\eqfinv
	\label{eq:new_power_flow}%
	\\
	&\additional_{i,a} = \sum_{j \in \NeighborsA{i}} \powerflow_{i,j}
	\eqsepv \forall a \in \Additional(i)
	\eqfinv
	\label{eq:new_variables}%
\end{align}
\end{subequations}
so that each constraint in~\eqref{eq:new_power_flow}-\eqref{eq:new_variables}
involves less variables than the original
aggregated formulation~\eqref{eq:power_flow}
--- assuming that $\card{\Additional(i)} < \card{\Neighbors{i}}$ .
Next, we add $2 \card{\Additional(i)}$
real variables to the POP~\eqref{eq:POP}
and extend $\na{\bijR, \bijIm}$
so that 
\begin{equation}
\begin{cases}
\bijR(a) = \Real{\additional_{i,a}} \ ,
\\
\bijIm(a) = \Imaginary{\additional_{i,a}} \ ,
\end{cases} \forall a \in \Additional(i)
\eqfinp
\end{equation}
Finally, we define minimal subsets
in the same spirit of~\eqref{eq:minimal_sparsity}:
\begin{subequations}
\label{eq:new_sparsity_total}
\begin{align}
\IndicesMS_{\#i} =& 
\na{\bijR(i), \bijIm(i)}
\bigcup_{a \in \Additional\np{i}} \na{\bijR(a), \bijIm(a)} \notag \\
&\bigcup_{g \in \Generators(i)} \na{\bijR(g), \bijIm(g)}
\eqfinv
\label{eq:new_sparsity_sets}%
\\
\IndicesMS_{\#a} =& 
\na{\bijR(i), \bijIm(i)}
\bigcup_{j \in \NeighborsA{i}} \na{\bijR(j), \bijIm(j)} \notag \\
&\bigcup \na{\bijR(a), \bijIm(a)}
\eqsepv \forall a \in \Additional(i)
\eqfinp
\label{eq:additional_sets}%
\end{align}
\end{subequations}

In turn, the sets  $\Additional\np{i}$
and  $\na{\NeighborsA{i}}_{a \in \Additional\np{i}}$
should be designed carefully to control the cardinalities
of the subsets defined 
by~\eqref{eq:new_sparsity_sets}-\eqref{eq:additional_sets}.
We suggest to use the solutions of the integer program
\begin{subequations}
\begin{equation}
	\label{eq:integer_program}%
	\min_{(n_\Additional, \overline{n}_a) \in \NN^{*2}} n_\Additional
	\quad \text{ s.t. } 
	\begin{cases}
		2(\overline{n}_a + 2) \leq \Imax \eqfinv
		\\
		n_\Additional \times \overline{n}_a
		\geq \card{\Neighbors{i}} \eqfinv
	\end{cases}
\end{equation}
which admits 
\begin{equation}
	n_\Additional = \left\lceil \frac{\card{\Neighbors{i}}}{\left\lfloor\frac{\Imax}{2}\right\rfloor - 2}
	\right\rceil
	\quad \text{ and } \quad
	\overline{n}_a = \left\lceil
	\frac{\card{\Neighbors{i}}}{n_\Additional}
	\right\rceil 
	\label{eq:integer_solution}
\end{equation}
\end{subequations}
as a solution, if $\Imax \geq 6$.
The rationale behind the formulation of Problem~\eqref{eq:integer_program}
is that we want 
to minimize $n_\Additional = \card{\Additional\np{i}}$
so as to reduce the cardinality of $\IndicesMS_{\#i}$
in~\eqref{eq:new_sparsity_sets}.
Meanwhile, we want to dispatch neighbors equally 
over the partition 
$\na{\NeighborsA{i}}_{a \in \Additional\np{i}}$,
which is composed of sets whose cardinalities are at most $\overline{n}_a$.
The constraints of Problem~\eqref{eq:integer_program}
ensure that the subsets 
$\na{\IndicesMS_{\#a}}_{a \in \Additional\np{i}}$
in~\eqref{eq:additional_sets} have cardinalities lower than
$\Imax$ (first inequality)
and that the partition
$\na{\NeighborsA{i}}_{a \in \Additional\np{i}}$
covers $\Neighbors{i}$ (second inequality).

Applying the solution~\eqref{eq:integer_solution},
we obtain a reduction of the cardinality of
$\IndicesMS_{\#i}$ provided that 
\begin{equation}
	\Imax \geq 4 + 
	\frac{2\card{\Neighbors{i}}}{\card{\Neighbors{i}} - 1}
	\eqfinp
\end{equation}
However, we might still have that 
$\card{\IndicesMS_{\#i}} > \Imax$. In this case,
we can operate a similar partitioning of $\Generators(i)$
to reduce the contribution of
power generation variables $\na{\power_g}_{g \in \Generators(i)}$
to the cardinality of $\IndicesMS_{\#i}$
in~\eqref{eq:new_sparsity_sets}.


\section{Numerical examples}
\label{sec:numerical_examples}

We illustrate the success of minimal sparsity in computing 
second-order moment relaxation bounds in AC-OPF.
In our experiments, we use Mosek 9.3 \cite{aps2020mosek}
to solve SDPs
and Ipopt \cite{wachter2006implementation}
for nonlinear programs. Both solvers are
applied with their default settings.
We display the results
of sparse SOS relaxations, i.e. the dual of \eqref{eq:sparse_moment_hierarchy},
as they are usually 
better handled than moment relaxations
by Mosek \cite[\S 7.5]{aps2020mosek}. 
The interface between data, models and solvers is implemented
with JuMP \cite{dunning2017jump} 
and PowerModels \cite{coffrin2018powermodels}.
We run experiments on a 2.10 GHz Intel CPU with 150 GB of RAM.
Our code is publicly available\footnote{\href{https://github.com/adrien-le-franc/MomentSOS.jl}{https://github.com/adrien-le-franc/MomentSOS.jl}} 
and we use open data from~\cite{godard2019resolution, babaeinejadsarookolaee2019power}.

We measure the accuracy of a relaxation in term of its
optimality gap
\begin{equation}
\gap_r(\Indices) = \frac{\valUb - \valSOS_r(\Indices)}{\valUb} \times 100 \eqfinv
\label{eq:gap}%
\end{equation}
where $\valUb$ is an upper bound computed with 
Ipopt. 
First, in \S\ref{sec:case57}, we measure the accuracy of
minimal sparsity on
modified case 57 IEEE instances 
that display large optimality gaps for first-order SOS relaxations.
Second, in \S\ref{sec:pglib}, we investigate
the scalability of minimal sparsity
on larger PGLib instances.

\subsection{Case 57 IEEE modified}
\label{sec:case57}

%

We consider the 1000 modified instances
generated in \cite[\S 5.4]{godard2019resolution}
by drawing random linear cost
parameters $\na{C_{1,g}}_{g \in \Generators}$
in~\eqref{eq:AC-OPF}
for case 57 IEEE from \cite{babaeinejadsarookolaee2019power}.
We concentrate on the ten instances
displaying the largest optimality gaps 
at the first-order SOS relaxation.
Following the formulation of~\cite{godard2019resolution},
we adopt a simplified AC-OPF model for this
experiment: 
limits on power lines and angle differences
in~\eqref{eq:line_constraints}-\eqref{eq:angle_constraints}
are ignored and we consider only linear costs.
Since moreover
case 57 IEEE has at most one generator per bus,
we may 
consider a voltage-only formulation
of Problem~\eqref{eq:AC-OPF}, and,
for the sake of numerical stability, 
we scale all polynomial coefficients to 
$f_{k,\alpha} \in \nc{-1, 1}$.
 
We present numerical results obtained with
clique-based sparsity in Table~\ref{tab:case57_clique}
and with minimal sparsity in Table~\ref{tab:case57_minimal}.
In both cases, we report the computing time
(columns 2-3) and the optimality gaps
(columns 4-5)
of the first- and second-order SOS relaxation bounds.

\paragraph{Bound accuracy}

we observe that the second-order relaxation
based on minimal sparsity always achieves zero
 optimality gap for all of the modified case 57 IEEE instances
(Table~\ref{tab:case57_minimal}, column 5).
This suggests that, despite its heuristic nature, minimal
sparsity is suitable to compute tight lower bounds for AC-OPF.
In turn, clique-based sparsity performs equally well
for the second-order relaxation 
(Table~\ref{tab:case57_clique}, column 5).
Interestingly, for first-order sparse relaxations,
the optimality gaps obtained with clique-based
sparsity (Table~\eqref{tab:case57_clique}, column 4) 
are smaller than the ones of minimal sparsity
(Table~\eqref{tab:case57_minimal}, column 4).
We also note that minimal sparsity yields more stable relaxations
than clique-based sparsity for instances 84 and 829,
for which the computation of $\valSOS_2(\IndicesCB)$
stopped with Mosek's \texttt{SLOW\_PROGRESS} termination status
(Table~\ref{tab:case57_clique}, column 5).
Better numerical stability could 
arise from the fact that 
minimal sparsity typically features smaller SDPs than 
clique-based sparsity.

\paragraph{Computing time}

the main improvement of minimal-sparsity
over a clique-based approach lies in the reduction
of computing time.
Indeed, evaluating clique-based second-order sparse relaxation bounds
$\valSOS_2(\Indices^\text{c})$ requires 3-6 hours
of computation per instance 
(Table~\ref{tab:case57_clique}, column 3),
whereas each of their minimal sparsity counterparts
$\valSOS_2(\IndicesMS)$
can be computed within one minute 
(Table~\ref{tab:case57_minimal}, column 3).
We believe that this shrinkage of computing time is due to the reduction of the size of the largest subsets in $\Indices$:
with clique-based sparsity, we have
$\max_p(\card{\IndicesCB_p}) = 26$,
while minimal sparsity features smaller cardinalities with
$\max_p(\card{\IndicesMS_p}) = 14$.
Lastly,
we mention that this way of certifying optimality
gaps
also outperforms
the branch-and-bound technique tested in~\cite{godard2019resolution},
which achieves an average of $0.16\%$ optimality gap
after 120 hours of computation per instance.

\begin{table}[!ht]
	\centering
	\begin{tabular}{c|r|r||r|r}
		\multirow{2}{*}{instances} & 
		\multicolumn{2}{c}{time (s)} &
		\multicolumn{2}{c}{gap (\%)} \\ 
		& $\valSOS_1(\Indices^\text{c})$  & $\valSOS_2(\Indices^\text{c})$
		& $\gap_1(\Indices^\text{c})$ & $\gap_2(\Indices^\text{c})$ \\ \hline
		84 & 6.88 $10^{\text{-}1}$ & 1.97 $10^{4}$ 
		& 3.05 & $^*$0.00 \\	
		260 & 5.90 $10^{\text{-}1}$ & 1.19 $10^{4}$ 
		& 1.67 & 0.00 \\
		267 & 6.54 $10^{\text{-}1}$ & 1.37 $10^{4}$ 
		& 1.21 & 0.00 \\
		299 & 6.17 $10^{\text{-}1}$ & 2.23 $10^{4}$ 
		& 1.92 & 0.00 \\
		391 & 5.53 $10^{\text{-}1}$ & 1.67 $10^{4}$ 
		& 1.25 & 0.00 \\
		628 & 6.16 $10^{\text{-}1}$ & 1.79 $10^{4}$ 
		& 6.64 & 0.00 \\
		683 & 6.80 $10^{\text{-}1}$ & 1.44 $10^{4}$ 
		& 2.32 & 0.00 \\
		829 & 6.39 $10^{\text{-}1}$ & 1.98 $10^{4}$ 
		& 2.00 & $^*$0.00 \\
		868 & 6.43 $10^{\text{-}1}$ & 1.41 $10^{4}$ 
		& 2.17 & 0.00 \\
		974 & 6.80 $10^{\text{-}1}$ & 1.45 $10^{4}$ 
		& 1.92 & 0.00
	\end{tabular}
	\caption{
		Results for AC-OPF case 57 modified
		and clique-based subsets ($\card{\IndicesCB} = 38 \eqsepv \max_p\card{\IndicesCB_p} = 26$).
		In all cases, the solution returned is primal feasible.
		Instances for which Mosek terminated
		with the \texttt{SLOW\_PROGRESS} status are marked
		with ``$*$''
	}
	\label{tab:case57_clique}
\end{table}

\begin{table}[!ht]
	\centering
	\begin{tabular}{c|r|r||r|r}
		\multirow{2}{*}{instances} & 
		\multicolumn{2}{c}{time (s)} &
		\multicolumn{2}{c}{gap (\%)} \\ 
		& $\valSOS_1(\Indices^\text{m})$  & $\valSOS_2(\Indices^\text{m})$
		& $\gap_1(\Indices^\text{m})$ & $\gap_2(\Indices^\text{m})$ \\ \hline
		84 & 1.91 $10^{\text{-}1}$ & 4.57 $10^{1}$ 
		& 3.30 & 0.00 \\
		260 & 2.58 $10^{\text{-}1}$	& 4.35 $10^{1}$ 
		&	1.85 & 0.00	\\
		267 & 1.86 $10^{\text{-}1}$	& 4.39 $10^{1}$ 
		&	1.42 & 0.00 \\
		299 & 1.96 $10^{\text{-}1}$ & 5.67 $10^{1}$ 
		&	2.06 & 0.00	\\
		391 & 1.83 $10^{\text{-}1}$	& 5.27	$10^{1}$ 
		&	1.54 & 0.00	\\
		628 & 2.00 $10^{\text{-}1}$ & 5.31	$10^{1}$ 
		&	6.89 & 0.00	\\
		683 & 1.95 $10^{\text{-}1}$ & 4.56	$10^{1}$ 
		&	2.50 & 0.00 \\
		829 & 1.84 $10^{\text{-}1}$	& 4.28	$10^{1}$ 
		&	2.21 & 0.00	\\
		868 & 1.92 $10^{\text{-}1}$	& 4.30 $10^{1}$ 
		&	2.33 & 0.00 \\
		974 & 2.08 $10^{\text{-}1}$	& 5.10 $10^{1}$ 
		&	2.08 & 0.00	
	\end{tabular}
	\caption{
		Results for AC-OPF case 57 modified
		and minimal sparsity
		subsets ($\card{\IndicesMS} = 57 \eqsepv \max_p(\card{\IndicesMS_p}) = 14$).
		In all cases, the solution returned is primal feasible
	}
	\label{tab:case57_minimal}
\end{table}

\subsection{Standard PGLib examples}
\label{sec:pglib}

We present further results on the standard
AC-OPF formulation~\eqref{eq:AC-OPF}.
To investigate on the scalability of the results obtained
on case 57 IEEE,
we consider all PGLib cases with up to 1000 buses
and large scale RTE cases with thousands of buses.
As larger instances tend to be less numerically stable,
we scale both polynomial coefficients to 
$f_{k,\alpha} \in \nc{-1, 1}$
and POP variables to $x_n \in \nc{0, 1}$.
We report the performance of 
second-order relaxations based on minimal sparsity
in Table~\ref{tab:minimal_sparsity_pglib}.
We apply a maximal cardinality threshold $\Imax$ = 12
as introduced in~\S\ref{sec:finer_control},
resulting in additional POP variables
(Table~\ref{tab:minimal_sparsity_pglib}, column 8).

\begin{table*}[!ht]
	\centering
	\begin{tabular}{l|r|r|r||r|r|r||r|r}
		\multirow{2}{*}{PGLib cases} & 
		\multicolumn{3}{c}{ computing time (s) } &
		\multicolumn{3}{c}{$\gap_2(\Indices^\text{m})$ gap (\%)} & 
		\multicolumn{2}{c}{POP variables} \\ 
		& TYP & API & SAD &
		TYP & API & SAD &   
		added & total \\ \hline		
		3 LMBD 		&	5.32 $10^{\text{-1}}$	&	9.78 $10^{\text{-1}}$	&	5.54 $10^{\text{-1}}$	&	0.00	&	0.00	&	0.00	& 0 & 12 \\
		5 PJM   	&	7.22 $10^0$	&	9.45 $10^0$	&	5.34 $10^0$	&	0.00				&	$^*$0.07	&	0.00	& 0 & 20 \\
		14 IEEE		&	2.03 $10^1$	&	2.52 $10^1$	&	1.89 $10^1$	&	0.00				&	0.00	&	0.00	& 0 & 38 \\
		24 IEEE RTS &	5.19 $10^1$	&	8.73 $10^1$	&	5.56 $10^1$	&	$^*$0.00			&	$^*$0.00	&	$^*$0.00	& 24 & 138 \\
		30 AS 		&	3.04 $10^1$	&	6.26 $10^1$	&	3.26 $10^1$	&	0.00				&	?	&	0.00	& 8 & 80 \\
		30 IEEE  	&	3.86 $10^1$	&	4.29 $10^1$	&	4.45 $10^1$	&	$^*$0.00			&	$^*$0.00	&	$^*$0.00	& 8 & 80 \\
		39 EPRI		&	2.45 $10^1$	&	2.36 $10^1$	&	2.25 $10^1$	&	0.00				&	$^*$0.16	&	0.00	& 0 & 98 \\
		57 IEEE		&	4.69 $10^1$	&	4.01 $10^1$	&	4.25 $10^1$	&	0.00				&	0.00	&	0.00	& 12 & 140 \\
		60 C 		&	5.94 $10^1$	&	1.21 $10^2$	&	8.92 $10^1$	&	$^*$0.00			&	$^*$0.09	&	$^*$0.03	& 20 & 186 \\
		73 IEEE RTS	&	1.85 $10^2$	&	2.76 $10^2$	&	1.93 $10^2$	&	$^*$0.00			&	$^*$\textbf{4.67}	&	$^*$0.05	& 80 & 424 \\
		89 PEGASE 	&	8.59 $10^2$	&	9.14 $10^2$	&	8.12 $10^2$	&	$^*$0.00			&	?	&	$^*$0.00	& 158 & 360 \\
		118 IEEE  	&	2.82 $10^2$	&	4.95 $10^2$	&	3.79 $10^2$	&	$^*$0.00			&	$^*$\textbf{9.59}	&	$^*$0.05	& 78 & 422 \\
		162 IEEE DTC&	9.51 $10^2$	&	1.21 $10^3$	&	1.09 $10^3$	&	$^*$0.70			&	$^*$0.44	&	$^*$0.41	& 88 & 436 \\
		179 GOC		&	2.61 $10^2$	&	3.63 $10^2$	&	2.83 $10^2$	&	$^*$0.04			&	$^*$0.35	&	$^*$\textbf{4.08}	& 100 & 516 \\
		200 ACTIV  	&	3.33 $10^2$	&	5.06 $10^2$	&	3.15 $10^2$	&	$^*$0.00			&	$^*$0.00	&	$^*$0.00	& 48 & 524 \\
		240 PSERC  	&	9.05 $10^2$	&	1.38 $10^3$	&	1.02 $10^3$	&	$^*$\textbf{1.61}	&	$^*$0.24	&	$^*$\textbf{2.94}	& 292 & 1058 \\
		300 IEEE  	&	9.35 $10^2$	&	8.66 $10^2$	&	1.01 $10^3$	&	$^*$0.00			&	$^*$0.29	&	$^*$0.00	& 68 & 806 \\
		500 GOC		&	1.59 $10^3$	&	1.92 $10^3$	&	1.45 $10^3$	&	$^*$0.00			&	$^*$\textbf{2.21}	&	$^*$\textbf{3.66}	& 308 & 1650 \\
		588 SDET  	&	1.08 $10^3$	&	8.50 $10^2$	&	1.15 $10^3$	&	$^*$0.25			&	?	&	$^*$0.22	& 106 & 1472 \\
		793 GOC  	&	1.09 $10^3$	&	1.22 $10^3$	&	1.31 $10^3$	&	?					&	?	&	$^*$\textbf{1.72}	& 116 & 1896 \\ \hline
		1888 RTE  	&	4.74 $10^3$	&	4.99 $10^3$	&	4.49 $10^3$	&	?					&	$^*$0.05	&	$^*$\textbf{2.58}	& 1048 & 5404 \\
		1951 RTE  	&	4.68 $10^3$	&	5.73 $10^3$	&	4.48 $10^3$	&	$^*$-0.01			&	$^*$0.18	&	$^*$0.25	& 1076 & 5710 \\
		2848 RTE  	&	6.07 $10^3$	&	8.38 $10^3$	&	6.48 $10^3$	&	?					&	?	&	?	& 1472 & 8190 \\
		2868 RTE  	&	6.98 $10^3$	&	7.62 $10^3$	&	6.85 $10^3$	&	?					&	?	&	$^*$0.39	& 1498 & 8356 \\
		6468 RTE  	&	1.27 $10^4$	&	1.92 $10^4$	&	1.49 $10^4$	&	$^*$0.27			&	?	&	?	& 3438 & 17172\\
		6470 RTE  	&	1.57 $10^4$	&	1.91 $10^4$	&	1.93 $10^4$	&	$^*$0.74			&	?	&	?	& 3478 & 17940\\
		6495 RTE  	&	1.53 $10^4$	&	1.71 $10^4$	&	1.71 $10^4$	&	?					&	?	&	?	& 3500 & 17850 \\
		6515 RTE 	&	1.81 $10^4$	&	1.88 $10^4$	&	1.53 $10^4$	&	?					&	?	&	?	& 3492 & 17890
	\end{tabular}
	\caption{
		Results on standard AC-OPF instances from PGLib
		for second-order relaxation gaps
		$\gap_2(\Indices^\text{m})$ based on minimal sparsity.
		The maximal cardinality threshold is set to $\Imax$ = 12.
		Instances for which Mosek
		terminated with the \texttt{SLOW\_PROGRESS} status
		are marked with ``$*$'' and ``$?$'' indicates
		an \texttt{UNKNOWN\_RESULT\_STATUS}}
	\label{tab:minimal_sparsity_pglib}
\end{table*}

\paragraph{Our results}

second-order minimal sparsity relaxations
successfully certify less than 1\%
of optimality gap for 47 of the 60 instances 
with up to 1000 buses
(Table~\ref{tab:minimal_sparsity_pglib}, columns 5-7).
For 8 other instances (with gap values in bold font),
Mosek stops at a feasible point
with the \texttt{SLOW\_PROGRESS} termination status,
which suggests that the accuracy of the bound 
$\valSOS_2(\Indices^\text{m})$ could be further improved.
Lastly, the solver returns an \texttt{UNKNOWN\_RESULT\_STATUS}
for the 5 other instances.

Addressing larger AC-OPF instances appears numerically
challenging, as Mosek stops with an \texttt{UNKNOWN\_RESULT\_STATUS}
for 16 out of the 24 large RTE instances
(Table~\ref{tab:minimal_sparsity_pglib}, columns 5-7). 
Moreover, we obtain a negative gap value for case 1951 RTE TYP,
which means that its bound $\valSOS_2(\Indices^\text{m})$
should be carefully certified.

Nevertheless, we manage to compute second-order relaxations
with less than 1\% of optimality gaps
for instances with over 6000 buses.
Remarkably, these bounds are computed within
the same computing time and memory resources
required to compute clique-based second-order bounds
for case 57 in Table~\ref{tab:case57_clique}.
We believe that these results are unprecedented
and encouraging,
as they open new perspectives for second-order relaxations
of large scale AC-OPF instances.

\paragraph{Comparison with other approaches}

we observe that minimal sparsity drastically reduces both computing
times and memory requirements compared with
the clique-based sparse second-order moment relaxations
reported in \cite[Table II, column time2]{gopinath2020proving}.
In the latter reference,
cases 89 PEGASE and 162 IEEE DTC return out-of-memory errors
on a computer allowed with 125 GB of RAM,
whereas our RAM usage peak
for all instances with no more than 1000 buses
is of 10 GB.

Due to numerical instabilities,
it is not straightforward to compare the tightness
of our optimality gaps with the ones
obtained with moment relaxations
in \cite{gopinath2020proving, wang2021certifying}
--- where in both situations, Mosek also terminates with 
\texttt{SLOW\_PROGRESS} or \texttt{UNKNOWN\_RESULT\_STATUS}
for many cases.

Still, there are 31 instances
for which both minimal sparsity and the 1.5 CS-TSSOS hierarchy
give reliable bounds in~\cite{wang2021certifying}
--- that is, instances 
for which Mosek returns a primal feasible solution.
For these instances, minimal sparsity gives a strictly smaller
(hence better)
optimality gap in 14 cases, 
and a strictly larger (hence worse) optimality 
gap in 8 cases.
As a concluding remark, we highlight that minimal sparsity
and CS-TSSOS need not be presented as competitors,
as minimal subsets $\IndicesMS$ could be advantageously used
to replace the clique-based ones involved
in the CS-TSSOS hierarchy \cite{wang2022cs}.

\section{Conclusion}

We have introduced minimal sparsity,
designed to improve the scalability of second-order
moment-SOS sparse relaxations of AC-OPF.
Our numerical test cases reveal that
minimal sparsity gives very accurate lower
bounds, while drastically reducing
the computing times and memory requirements
over standard clique-based sparse relaxations.
Our best achievement is to compute
second-order relaxation bounds
certifying less than 1\% of optimality gaps
for instances with over 6000 buses.
Yet, such large instances remain numerically challenging
for state-of-the-art SDP solvers
--- in line with the conclusions of \cite{wang2021certifying}.
Regarding future improvements,
we look forward to ongoing progresses
in SDP solvers, and 
pre- or post-processing techniques
enforcing numerical stability and robustness
of SDP relaxations, as presented, e.g., in 
\cite{oustry2022certified,mai2022exploiting}.

\section*{Acknowledgments}
\noindent 
This work was supported by the PEPS2 FastOPF funded by RTE and the French Agency for mathematics in interaction with industry and society (AMIES), the EPOQCS grant funded by the LabEx CIMI (ANR-11-LABX-0040), the European Union’s Horizon 2020 research and innovation programme under the Marie Sk{\l}odowska-Curie Actions, grant agreement 813211 (POEMA), by the AI Interdisciplinary Institute ANITI funding, through the French ``Investing for the Future PIA3'' program under the Grant agreement n${}^\circ$ ANR-19-PI3A-0004 as well as by the National Research Foundation, Prime Minister’s Office, Singapore under its Campus for Research Excellence and Technological Enterprise (CREATE) programme.

\bibliographystyle{IEEEtran}
\bibliography{biblio}

\begin{thebibliography}{10}
\providecommand{\url}[1]{#1}
\csname url@samestyle\endcsname
\providecommand{\newblock}{\relax}
\providecommand{\bibinfo}[2]{#2}
\providecommand{\BIBentrySTDinterwordspacing}{\spaceskip=0pt\relax}
\providecommand{\BIBentryALTinterwordstretchfactor}{4}
\providecommand{\BIBentryALTinterwordspacing}{\spaceskip=\fontdimen2\font plus
\BIBentryALTinterwordstretchfactor\fontdimen3\font minus
  \fontdimen4\font\relax}
\providecommand{\BIBforeignlanguage}[2]{{%
\expandafter\ifx\csname l@#1\endcsname\relax
\typeout{** WARNING: IEEEtran.bst: No hyphenation pattern has been}%
\typeout{** loaded for the language `#1'. Using the pattern for}%
\typeout{** the default language instead.}%
\else
\language=\csname l@#1\endcsname
\fi
#2}}
\providecommand{\BIBdecl}{\relax}
\BIBdecl

\bibitem{babaeinejadsarookolaee2019power}
S.~Babaeinejadsarookolaee, A.~Birchfield, R.~D. Christie, C.~Coffrin,
  C.~DeMarco, R.~Diao, M.~Ferris, S.~Fliscounakis, S.~Greene, R.~Huang
  \emph{et~al.}, ``{The Power Grid Library for Benchmarking AC Optimal Power
  Flow Algorithms},'' \emph{arXiv preprint arXiv:1908.02788}, 2019.

\bibitem{bukhsh2013local}
W.~A. Bukhsh, A.~Grothey, K.~I. McKinnon, and P.~A. Trodden, ``Local solutions
  of the optimal power flow problem,'' \emph{IEEE Transactions on Power
  Systems}, vol.~28, no.~4, pp. 4780--4788, 2013.

\bibitem{molzahn2019survey}
D.~K. Molzahn, I.~A. Hiskens \emph{et~al.}, ``A survey of relaxations and
  approximations of the power flow equations,'' \emph{Foundations and
  Trends{\textregistered} in Electric Energy Systems}, vol.~4, no. 1-2, pp.
  1--221, 2019.

\bibitem{josz2014application}
C.~Josz, J.~Maeght, P.~Panciatici, and J.-C. Gilbert, ``{Application of the
  moment-SOS approach to global optimization of the OPF problem},'' \emph{IEEE
  Transactions on Power Systems}, vol.~30, no.~1, pp. 463--470, 2014.

\bibitem{molzahn2014moment}
D.~K. Molzahn and I.~A. Hiskens, ``Moment-based relaxation of the optimal power
  flow problem,'' in \emph{2014 Power Systems Computation Conference}.\hskip
  1em plus 0.5em minus 0.4em\relax IEEE, 2014, pp. 1--7.

\bibitem{lasserre2001global}
J.-B. Lasserre, ``Global optimization with polynomials and the problem of
  moments,'' \emph{SIAM Journal on Optimization}, vol.~11, no.~3, pp. 796--817,
  2001.

\bibitem{nie2014optimality}
J.~Nie, ``{Optimality conditions and finite convergence of Lasserre’s
  hierarchy},'' \emph{Mathematical Programming}, vol. 146, pp. 97--121, 2014.

\bibitem{gopinath2020proving}
S.~Gopinath, H.~L. Hijazi, T.~Weisser, H.~Nagarajan, M.~Yetkin, K.~Sundar, and
  R.~W. Bent, ``Proving global optimality of acopf solutions,'' \emph{Electric
  Power Systems Research}, vol. 189, p. 106688, 2020.

\bibitem{waki2006sums}
H.~Waki, S.~Kim, M.~Kojima, and M.~Muramatsu, ``Sums of squares and
  semidefinite program relaxations for polynomial optimization problems with
  structured sparsity,'' \emph{SIAM Journal on Optimization}, vol.~17, no.~1,
  pp. 218--242, 2006.

\bibitem{molzahn2014sparsity}
D.~K. Molzahn and I.~A. Hiskens, ``Sparsity-exploiting moment-based relaxations
  of the optimal power flow problem,'' \emph{IEEE Transactions on Power
  Systems}, vol.~30, no.~6, pp. 3168--3180, 2014.

\bibitem{sparsebook}
V.~Magron and J.~Wang, ``Sparse polynomial optimization: theory and practice,''
  \emph{Series on Optimization and Its Applications, World Scientific Press},
  2023, to appear.

\bibitem{wang2022cs}
J.~Wang, V.~Magron, J.~B. Lasserre, and N.~H.~A. Mai, ``{CS-TSSOS: Correlative
  and term sparsity for large-scale polynomial optimization},'' \emph{ACM
  Transactions on Mathematical Software}, vol.~48, no.~4, pp. 1--26, 2022.

\bibitem{wang2021certifying}
J.~Wang, V.~Magron, and J.~B. Lasserre, ``{Certifying global optimality of
  AC-OPF solutions via sparse polynomial optimization},'' \emph{Electric Power
  Systems Research}, vol. 213, p. 108683, 2022.

\bibitem{lasserre2006convergent}
J.-B. Lasserre, ``{Convergent SDP-relaxations in polynomial optimization with
  sparsity},'' \emph{SIAM Journal on Optimization}, vol.~17, no.~3, pp.
  822--843, 2006.

\bibitem{nesterov1994interior}
Y.~Nesterov and A.~Nemirovskii, \emph{Interior-point polynomial algorithms in
  convex programming}.\hskip 1em plus 0.5em minus 0.4em\relax SIAM, 1994.

\bibitem{bienstock2019strong}
D.~Bienstock and A.~Verma, ``{Strong NP-hardness of AC power flows
  feasibility},'' \emph{Operations Research Letters}, vol.~47, no.~6, pp.
  494--501, 2019.

\bibitem{bienstock2020mathematical}
D.~Bienstock, M.~Escobar, C.~Gentile, and L.~Liberti, ``Mathematical
  programming formulations for the alternating current optimal power flow
  problem,'' \emph{4OR}, vol.~18, no.~3, pp. 249--292, 2020.

\bibitem{molzahn2013implementation}
D.~K. Molzahn, J.~T. Holzer, B.~C. Lesieutre, and C.~L. DeMarco,
  ``Implementation of a large-scale optimal power flow solver based on
  semidefinite programming,'' \emph{IEEE Transactions on Power Systems},
  vol.~28, no.~4, pp. 3987--3998, 2013.

\bibitem{sliwak2020clique}
J.~Sliwak, E.~D. Andersen, M.~F. Anjos, L.~L{\'e}tocart, and E.~Traversi, ``A
  clique merging algorithm to solve semidefinite relaxations of optimal power
  flow problems,'' \emph{IEEE Transactions on Power Systems}, vol.~36, no.~2,
  pp. 1641--1644, 2020.

\bibitem{sliwak2021resolution}
J.~Sliwak, ``R{\'e}solution de probl{\`e}mes d'optimisation pour les
  r{\'e}seaux de transport d'{\'e}lectricit{\'e} de grande taille avec des
  m{\'e}thodes de programmation semi-d{\'e}finie positive,'' Ph.D.
  dissertation, Polytechnique Montr{\'e}al, 2021.

\bibitem{magron2021tssos}
\BIBentryALTinterwordspacing
V.~Magron and J.~Wang, ``{TSSOS: a Julia library to exploit sparsity for
  large-scale polynomial optimization},'' \emph{Proceedings of {\it MEGA:
  Effective Methods in Algebraic Geometry}}, 2021. [Online]. Available:
  \url{https://github.com/wangjie212/TSSOS}
\BIBentrySTDinterwordspacing

\bibitem{bodlaender2010treewidth}
H.~L. Bodlaender and A.~M. Koster, ``{Treewidth computations I. Upper
  bounds},'' \emph{Information and Computation}, vol. 208, no.~3, pp. 259--275,
  2010.

\bibitem{aps2020mosek}
M.~ApS, ``Mosek modeling cookbook,'' 2020.

\bibitem{wachter2006implementation}
A.~W{\"a}chter and L.~T. Biegler, ``On the implementation of an interior-point
  filter line-search algorithm for large-scale nonlinear programming,''
  \emph{Mathematical Programming}, vol. 106, pp. 25--57, 2006.

\bibitem{dunning2017jump}
I.~Dunning, J.~Huchette, and M.~Lubin, ``{JuMP: A modeling language for
  mathematical optimization},'' \emph{SIAM review}, vol.~59, no.~2, pp.
  295--320, 2017.

\bibitem{coffrin2018powermodels}
C.~Coffrin, R.~Bent, K.~Sundar, Y.~Ng, and M.~Lubin, ``Powermodels. jl: An
  open-source framework for exploring power flow formulations,'' in \emph{2018
  Power Systems Computation Conference (PSCC)}.\hskip 1em plus 0.5em minus
  0.4em\relax IEEE, 2018, pp. 1--8.

\bibitem{godard2019resolution}
H.~Godard, ``R{\'e}solution exacte du probl{\`e}me de l'optimisation des flux
  de puissance,'' Ph.D. dissertation, Paris, CNAM, 2019.

\bibitem{oustry2022certified}
A.~Oustry, C.~D’Ambrosio, L.~Liberti, and M.~Ruiz, ``{Certified and accurate
  SDP bounds for the ACOPF problem},'' \emph{Electric Power Systems Research},
  vol. 212, p. 108278, 2022.

\bibitem{mai2022exploiting}
N.~H.~A. Mai, J.-B. Lasserre, V.~Magron, and J.~Wang, ``Exploiting constant
  trace property in large-scale polynomial optimization,'' \emph{ACM
  Transactions on Mathematical Software}, vol.~48, no.~4, pp. 1--39, 2022.

\end{thebibliography}

\end{document}